\documentclass[12pt]{amsart}
\usepackage{amssymb}
\usepackage{amscd}
\usepackage[a4paper]{geometry}
\usepackage{graphics}

\usepackage[
hyperref,
]{degt}

\def\CG#1{\Z/#1}
\def\gcd{\QOPNAME{gcd}}
\def\lcm{\QOPNAME{lcm}}
\def\units{^{\times\!}}
\def\fA{\frak{A}}
\def\fB{\frak{B}}

\let\Fermat=\Phi

\let\rl=R

\let\ln=L

\let\dv=V
\let\e\epsilon

\def\two{_=}

\def\sset{\bold S}
\def\tset{\bold T}
\def\Kset{\bold K}
\def\Tset{\bold T}
\def\ofD#1{\<#1\>}
\def\SS{\sset\ofD}
\def\KK{\bold K\ofD}
\def\TT{\tset\ofD}
\def\QQ{\bold Q}
\def\GG{\Bbb G}
\def\Pic{\QOPNAME{Pic}}

\def\Log{\QOPNAME{Log}}
\def\height{\QOPNAME{ht}}

\def\SF#1{\sset\of{#1}}
\def\KF#1{\Kset\of{#1}}
\def\TF#1{\Tset\of{#1}}
\def\FF#1{\Fermat\of{#1}}
\def\FC#1{\Fermat^\circ\of{#1}}

\def\CJ{\Cal J}
\def\CK{\Cal K}
\def\CS{\Cal S}

\def\of#1{(#1)}
\def\of#1{[#1]}

\def\tX{\smash{\tilde X}}

\def\Am{\mathrm{A}\of}

\def\Bm{\mathrm{B}\of}

\def\Hc{H_\circ}
\def\dual{^\vee}

\let\=\B

\let\vg=v
\let\hg=h
\let\ag=a
\let\cg=c
\let\ug=u
\let\g=g
\let\3\tilde

\let\GL\Lambda

\let\inj\inc

\def\CW{{\sl CW}}

\title{Projective spaces in Fermat varieties}

\author{Alex Degtyarev}

\address{%
Department of Mathematics\\
Bilkent University\\
06800 Ankara, Turkey}

\email{degt@fen.bilkent.edu.tr}

\thanks{%
The author was partially supported by the JSPS grant L15517
and T\"{U}B\DOTaccent{I}TAK grant 114F325}



\begin{document}

\begin{abstract}
We give a brief systematic overview of a few
results concerning the
N\'eron--Severi lattices of Fermat varieties and
Delsarte surfaces.
\end{abstract}

\maketitle

\tableofcontents

\section{Introduction}\label{S.intro}

The goal of this survey is to give a brief systematic overview of a few
results, both recent and old, concerning the generators of the
N\'eron--Severi lattices of Fermat varieties and of closely related to them
Delsarte surfaces.

Citing T.~Shioda~\cite{Shioda:Fermat1},
the N\'{e}ron--Severi group
``\dots is a rather delicate invariant of arithmetic nature.
Perhaps for this reason
it usually requires some nontrivial work before one can determine the Picard
number of a given variety, let alone the full structure of its N\'{e}ron--Severi
group.''
The Picard ranks $\rank\NS(X)$
of Fermat varieties and Delsarte surfaces were computed
in~\cite{Shioda:Fermat,Shioda:var1,Shioda:Delsarte}; these results are
outlined in \autoref{s.Fermat}. Comparing $\rank\NS(X)$ with the rank of the
subgroup $\sset\subset\NS(X)$ generated by a certain set $\CS$ of
``immediately seen''
subvarieties of~$X$ (projective
spaces or images thereof, see \autoref{s.planes} and \autoref{s.Delsarte}),
it was observed that, under some rather general assumptions, the subvarieties
constituting~$\CS$ generate
the rational N\'eron--Severi group: $\sset\otimes\Q=\NS(X;\Q)$.
Naturally, the question arose whether one also has $\sset=\NS(X)$ over the
integers; the affirmative answer to this question would give one the complete
structure of the N\'eron--Severi lattice.

The question remained unsettled for almost 30 years, until the first numerical
evidence suggesting the positive answer appeared in 2010,
see~\cite{Shioda:Fermat2,Shimada.Takahashi:primitivity}.
The original case of Fermat surfaces was finally settled (in the affirmative)
in~\cite{degt:Fermat}. The situation with Delsarte surfaces turned out more
complicated: it was shown in~\cite{degt:Delsarte} that the answer depends on
the structure of the defining equation and typically is in the negative,
although the torsion of the quotient $\NS(X)/\sset$ is bounded; \eg, its
length does not exceed~$7$.
The techniques used in the proofs are outlined in
\autoref{s.torsion} and \autoref{s.AM}, and a brief
account of the results is found in \autoref{s.partial}.

The most recent achievement is an algebraic restatement
(similar to that used in \cite{degt:Fermat}) of the original
question for Fermat varieties of
higher dimension, see~\cite{degt.Shimada} and \autoref{s.reduction}: the
answer is given in terms of the integral torsion of certain modules over
polynomial rings. Unfortunately, we failed to prove that this torsion
vanishes. So far, only some numerical evidence and a few partial vanishing
results are available, see \autoref{s.vanishing}. Some of these
partial results have
geometric implications to a wider class of varieties; they are discussed in
\autoref{s.others}.

In \autoref{S.problems}, I briefly state a few open problems
that seem to be of general interest.

\subsection*{Acknowledgements}
I would like to thank the organizers of the Kinosaki Algebraic Geometry
Symposium who gave me the opportunity to present my work at this prestigious
forum. I am grateful to my co-author Ichiro Shimada, who has always
generously shared
his knowledge and ideas. My special gratitude goes to Tetsuji Shioda for his
continued support and motivation and a number of fruitful and instructive
discussions. This paper was written during my stay at Hiroshima
University, supported by the Japan Society for the Promotion of Science.

\section{The rational N\'eron--Severi lattice}\label{S.NSQ}

\subsection{Fermat varieties}\label{s.Fermat}
Consider the Fermat variety~$\Fermat_m^d\subset\Cp{d+1}$ given by
\[*
z_0^m+\ldots+z_{d+1}^m=0,
\]
and let $G_m$ be the group
\[
G_m:=\bigl\{(\e_0,\ldots,\e_{d+1})\in(\C\units)^{d+2}\bigm|
 \e_0^m=\ldots=\e_{d+1}^m=1\bigr\}/\text{diagonal}.
\label{eq.Gm}
\]
This group acts on~$\Fermat_m^d$ \via\
\[*
(\e_0,\ldots,\e_{d+1})\:(z_0:\ldots:z_{d+1})\mapsto(\e_0z_0:\ldots:\e_{d+1}z_{d+1}),
\]
inducing a decomposition
$H_2(\Fermat_m^d;\C)=\bigoplus_\Go H_d^\Go(\Fermat_m^d)$, where $\Go$ runs
through the dual group $G_m^\vee:=\Hom(G_m,\C\units)$,
\[*
G_m^\vee=\bigl\{(\Go_0,\ldots,\Go_{d+1})\in(\C\units)^{d+2}\bigm|
 \Go_0^m=\ldots=\Go_{d+1}^m=\Go_0\ldots\Go_{d+1}=1\bigr\}.
\]
According to~\cite{Shioda:var1}, one has $H_d^\Go=0$ unless either $\Go=1$ or
\[*
\Go\in\fA_m^d:=\bigl\{\Go\in G_m^\vee\bigm|
 \text{$\Go_i\ne1$ for each $i=0,\ldots,d+1$}\bigr\}.
\]
Furthermore, the dimension of each nontrivial eigenspace is~$1$
and the Hodge weight of the eigenspace $H_d^\Go$ corresponding to
a character $\Go\in\fA_m^d$ equals
$\Log\Go-1$, where
$\Log\Go:=\sum_i\Log\Go_i\in\{1,\ldots,d\}$
and $\Log\Go_i$ is the argument of the complex number~$\Go_i$
\emph{specialized to $[0,2\pi)$} and divided by~$2\pi$.
(Note a mysterious similarity between these formulas and the signature of a
generalized Hopf link, see~\cite{DFL}; I do not know a conceptual explanation
of this fact.)

The group of units $(\CG{m})\units$ acts on the character group $G_m^\vee$
\via\ $u\:(\Go_i)\mapsto(\Go_i^u)$
and it is clear that a sum $\bigoplus_{\Go}H_d^\Go$,
$\Go\in\Omega$,
of eigenspaces
is defined over~$\Q$ if and only if the index set $\Omega\subset G_m^\vee$
is invariant under this
action.
Therefore, assuming that $d=2k$ is even, the dimension of
the rational N\'eron--Severi lattice
\[*
\NS(\Fermat_m^d;\Q):=\NS(\Fermat_m^d)\otimes\Q
 =H_d(\Fermat_m^d;\Q)\cap H^{k,k}(\Fermat_m^d)
\]
equals $\ls|\fB_m^d|+1$, where
\[*
\fB_m^d:=\bigl\{\Go\in\fA_m^d\bigm|
 \text{$\Log(\Go^u)=k+1$ for each $u\in(\CG{m})\units$}\bigr\},
\]
see \cite{Shioda:var1,Shioda:Fermat}.

In the case of surfaces ($d=2$), the set $\fB_m^2$ has been studied in
\cite{Shioda:Fermat}. In particular, it has been shown that
\begin{multline}
\quad\dim\NS(\Fermat_m^2;\Q)=3(m-1)(m-2)+\delta_m+1\\
 +24(m/3)^*+48(m/2)^*+24\epsilon(m),\qquad
\label{eq.dim.NS}
\end{multline}
where $\delta_m:=1-(m\bmod2)\in\{0,1\}$,
the expression $(q)^*$ stands for $q$ if $q\in\Z$ and $0$ otherwise, and
$\epsilon(m)$ is a bounded function that can be expressed as a certain sum
over the divisors $d\mathbin|m$ such that $\gcd(d,6)>1$ and $d\le180$.
Note that the last three terms vanish whenever $\gcd(m,6)=1$.


\subsection{Counting projective spaces}\label{s.planes}
Assume that $d=2k$ is even and
pick an unordered partition
\[
J=\bigl\{\{p_0,q_0\},\ldots,\{p_k,q_k\}\bigr\}
\label{eq.partition}
\]
of the index set $\{0,\ldots,d+1\}$ into $(k+1)$ unordered pairs. Then,
for each sequence $\eta=(\eta_0,\ldots,\eta_k)$
of $m$-th roots of $(-1)$, the projective $d$-space
\[
L_{J,\eta}:=\{z_{p_i}=\eta_iz_{q_i},\ i=0,\ldots,k\}
\label{eq.L}
\]
lies in $\Fermat_m^d$. Varying~$J$ and~$\eta$, we obtain $(2k+1)!!\,m^{k+1}$
distinct subspaces; their classes generate a certain subgroup
$\sset_m^d\subset\NS(X)$.

If $d=2$, the above spaces are lines and, for $m\ge3$, it can easily be shown
that there are no other lines in $\Fermat_m^2$ (see. \eg,
\cite{Boissiere.Sarti}). In this special case,
analyzing the intersection matrix (see, \eg, \cite{Shioda:Fermat2}), one
can also show that
\[
\rank\sset_m^2=3(m-1)(m-2)+\delta_m+1.
\label{eq.rank.S}
\]
(Alternatively, this rank can be found using \autoref{th.bounds} below.)
Comparing this to~\eqref{eq.dim.NS}, we arrive at the following statement.

\theorem[see~\cite{Shioda:Fermat}]\label{th.Shioda}
If $m\le5$ or $\gcd(m,6)=1$, then $\sset_m^2\otimes\Q=\NS(\Fermat_m^2;\Q)$.
\endtheorem

If $d\ge4$, a similar statement can be obtained by other means (induction
rather than direct counting).

\theorem[see~\cite{Ran,Shioda:var1}]\label{th.Shioda.var}
If $m=4$ or $m$ is prime, then $\sset_m^d\otimes\Q=\NS(\Fermat_m^d;\Q)$.
\endtheorem

Hence, a natural question, first raised in~\cite{Shioda:Fermat1}, is whether, under
the assumptions of Theorem~\ref{th.Shioda} or \ref{th.Shioda.var}, we have
the equality
$\sset_m^2=\NS(\Fermat_m^d)$, \ie, whether the
classes of the projective subspaces contained in~$\Fermat_m^d$ generate
$\NS(\Fermat_m^d)$ over the integers given that they do so over the
rationals. To answer this question, we will study the torsion group
$\tset_m^d:=\Tors\bigl(H_d(\Fermat_m^d)/\sset_m^d\bigr)$.
(Throughout the paper, the notation $\Tors A$ always stands for the integral
torsion of the abelian group~$A$, even if the latter is also a module over
another ring.)

\subsection{Delsarte surfaces}\label{s.Delsarte}
A \emph{Delsarte surface} $\Fermat_A\subset\Cp3$ is a surface given
by a four-term equation of the form
\[
\sum_{i=0}^3\prod_{j=0}^3z_j^{a_{ij}}=0,
\label{eq.Delsarte}
\]
see~\cite{Delsarte,Shioda:Delsarte},
where the exponent matrix~$A:=[a_{ij}]$ satisfies the following conditions:
\roster
\item\label{A.integer}
each entry~$a_{ij}$, $0\le i,j\le3$, is a non-negative integer;
\item\label{A.0}
each column of~$A$ has at least one zero;
\item\label{A.const}
$(1,1,1,1)$ is an eigenvector of~$A$,
\ie, $\sum_{j=0}^3a_{ij}=\Gl=\const(i)$;
\item\label{A.det}
$A$ is non-degenerate, \ie, $\det A\ne0$.
\endroster
Condition~\iref{A.0} asserts that the surface does not contain a coordinate
plane, and \iref{A.const} makes~\eqref{eq.Delsarte} homogeneous, the degree
being the eigenvalue~$\Gl$.

In general, this surface is singular and we silently replace~$\Fermat_A$ with
its resolution of singularities. The particular choice of the resolution is not
important as we will only deal with birational invariants.

Following~\cite{Shioda:Delsarte}, introduce the cofactor matrix
$A^*:=(\det A)A\1$ and let
\[*
d:=\gcd(a_{ij}^*),\quad
m:=\ls|\det A|/d,\quad
B=[b_{ij}]:=mA\1=\pm d\1A^*.
\]
Then, we have maps
\[*
\Fermat_m^2\overset{\pi_B}\longto\Fermat_A
\overset{\pi_A}\longto\Fermat:=\Fermat_1^2
\]
given by
\[*
\pi_B\:(z_i)\mapsto\Biggl(\prod_{j=0}^3z_j^{b_{ij}}\Biggr),\quad
\pi_A\:(z_i)\mapsto\Biggl(\prod_{j=0}^3z_j^{a_{ij}}\Biggr).
\]
Both maps are ramified coverings;
$\pi_A$ and $\pi_B\circ\pi_A\:(z_i)\mapsto(z_i^m)$
are ramified over the union
$\rl:=\rl_0+\rl_1+\rl_2+\rl_3\subset\Fermat$ of the
traces of the coordinate
planes, $\rl_i:=\Fermat\cap\{z_i=0\}$.
The $3m^2$ lines in~$\Fermat_m^2$ (see \autoref{s.planes}) project to the
three lines
\[*
L_{ij}:=\Fermat\cap\{z_i=z_j\},\quad 0\le i<j\le3.
\]
(Obviously, $L_{ij}=L_{kl}$ whenever $i,j,k,l$ are pairwise distinct, \ie,
the $L$-lines are indexed by partitions~$J$ as in \eqref{eq.partition}.)
Together, $R$ and $L:=L_{01}+L_{02}+L_{03}$ form the so-called \emph{Ceva-7}
arrangement in the projective plane~$\Fermat$, see \autoref{fig.divisor}
(where $R_0$ is the missing line at infinity).
\figure
\centerline{\cpic{lines}}
\caption{The divisor $\dv:=\ln+\rl\subset\Fermat$}\label{fig.divisor}
\endfigure

Since $R$ is a nodal curve,
the fundamental group $\GG:=\pi_1(\Fermat\sminus R)$ is abelian: it has four
generators $t_i$ dual to $[R_i]$, $i=0,\ldots,3$, that are
subject to the relation $t_0t_1t_2t_3=1$. The
finite ramified coverings $\pi_A$ as
above are in a natural one-to-one correspondence with \emph{finite quotients}
of $\GG$, \ie, epimorphisms $\Ga\:\GG\onto G$ onto a finite group~$G$.
Henceforth, we can disregard the original matrix~$A$ and speak about the
Delsarte surface $\FF\Ga$, which is defined as (any) smooth analytic
compactification of the covering of $\Fermat\sminus R$ corresponding to~$\Ga$.
In this notation,
$\Fermat_m^2=\FF{m}$, where an integer
$m\in\Z$ is regarded as a map $m\:\GG\onto\GG/m\GG$.

Found in~\cite{Shioda:Delsarte} is an algorithm making use
of~\eqref{eq.dim.NS} and computing the Picard rank (or rather corank, which is
a birational invariant)
of $\FF\Ga$ in terms
of~$\Ga$. On the other hand,
there is an ``obvious''
divisor $\dv[\Ga]:=\pi_A^*(R+L)$
in $\FF\Ga$ that plays the r\^ole of lines in~$\Fermat_m^2$;
let $\SF\Ga\subset\NS(\FF\Ga)$ be the subgroup generated by the
components of $\dv[\Ga]$.
One of the outcomes of~\cite{Shioda:Delsarte} is the equality
$\SF\Ga\otimes\Q=\NS(\FF\Ga;\Q)$ that holds whenever $\gcd(m,6)=1$
and the natural question whether, under the same assumption, one also has
$\SF\Ga=\NS(\FF\Ga)$ over the integers, \ie, whether the group
\[*
\TF\Ga:=\Tors\bigl(\NS(\FF\Ga)/\SF\Ga\bigr)
\]
is trivial.
Note that, in the case of a Fermat surface, this question is equivalent to
the original one raised in~\cite{Shioda:Fermat1}, see the end of
\autoref{s.Fermat}. Indeed, in this case, each divisorial pull-back
$\pi_A^*R_i$, $i=0,\ldots,3$, is a reduced
irreducible Fermat curve and,
fixing~$J$ and $\eta_1,\eta_2,\eta_3$ in~\eqref{eq.L}, we obtain $m$ lines
whose classes sum up to $[\pi_A^*R_i]$. Hence,
whenever $\Ga=m\in\N_+$, we have $\SF{m}=\sset_m^2$ and $\TF{m}=\tset_m^2$.

\section{The topological reduction}\label{S.reduction}

\subsection{The torsion group}\label{s.torsion}
Given a divisor~$D$ in a smooth compact surface~$X$, let
$\SS{D}\subset\NS(X)$ be the subgroup generated by the irreducible
components of~$D$.
Here and below, the N\'eron--Severi lattice $\NS(X)$ is the image of $\Pic X$
in the free abelian group $H_2(X)/\Tors$, which is canonically identified
with
$H^2(X)/\Tors$ \via\ Poincar\'e duality. The homomorphism is given by
$D\mapsto[D]$ in the language of divisors and homology or by
$\Cal L\mapsto c_1(\Cal L)$ in the language of line bundles and cohomology.
Thus,
\[*
\SS{D}:=\Im\bigl[\inj_*\:H_2(D)\to H_2(X)/\Tors\bigr],
\]
where $\inj\:D\into X$ is the inclusion.
We will also consider the groups
\[*
\KK{D}:=\Ker\bigl[\inj_*\:H_2(D)\to H_2(X)\bigr],\qquad
\TT{D}:=\Tors\bigl(\NS(X)/\SS{D}\bigr).
\]
The following statement is essentially the definition of $\Ext$ and
Poincar\'e duality.

\theorem[see \cite{degt:Delsarte,degt:Fermat}]\label{th.reduction}
For $D\subset X$ as above, let
\[*
K(X,D):=\Ker\bigl[\kappa_*\:H_1(X\sminus D)\to H_1(X)\bigr]
\]
be the kernel of the homomorphism~$\kappa_*$ induced by the inclusion
$\kappa\:X\sminus D\into X$.
Then there are
canonical
isomorphisms
\[*
\Tors K(X,D)=\Hom(\TT{D},\Q/\Z),\qquad
K(X,D)/\!\Tors=\Hom(\KK{D},\Z).
\]
\endtheorem

Indeed, $\kappa_*$
is Poincar\'{e} dual to
the homomorphism~$\rel$ in the exact sequence
\[*
\longto H^2(X)\overset{\inj^*}\longto H^2(D)\longto
H^3(X,D)\overset\rel\longto H^3(X)\longto.
\]
Thus, $K(X,D)=\Coker\inj^*$. The abelian group $H^2(D)$ is free and,
modulo torsion in $H^2(X)$, the homomorphism $\inj^*$ is the adjoint of
$\inj_*$ in the free resolution
\[*
\longto H_2(D)\overset{\inj_*}\longto H_2(X)/\Tors\longto T\longto0,
\]
where $T:=H_2(X)/\bigl(\Tors H_2(X)+\SS{D}\bigr)$. Thus,
\[*
\Coker\inj^*=\Ext(T,\Z)=\Hom(\Tors T,\Q/\Z).
\]
On the other hand, the quotient
$H_2(X)/c_1(\Pic X)$ is known to be torsion
free and we have $\Tors T=\TT{D}$.

Now, let $X=\FF\Ga$ be a Delsarte surface and $D=\dv[\Ga]$. To avoid
excessive nested parentheses, we will abbreviate
\[*
\SF\Ga:=\SS{\dv[\Ga]},\quad
\TF\Ga:=\TT{\dv[\Ga]},\quad
\KF\Ga:=\KK{\dv[\Ga]}
\]
and use the shortcut $\FC\Ga:=\FF\Ga\sminus\dv[\Ga]$.
The group $H_1(\FF\Ga)=\pi_1(\FF\Ga)$ is finite
abelian (see~\cite{degt:Fermat} or~\eqref{eq.pi1} below)
and the homomorphism~$\kappa_*$
in \autoref{th.reduction} factors
through the free abelian group
\[*
H_1(\FF\Ga\sminus\pi_A\1R)=\pi_1(\FF\Ga\sminus\pi_A\1R)=\Ker\Ga\subset\GG.
\]
Hence, \autoref{th.reduction} can be recast in a simpler form
\[
\Tors H_1(\FC\Ga)=\Hom(\TF\Ga,\Q/\Z),\qquad
H_1(\FC\Ga)/\Tors=\Hom(\KF\Ga,\Z).
\label{eq.torsion=H1}
\]
Note also that, as long as the torsion is concerned, $H_1(\FC\Ga)$ can be
replaced with $C_1/\Im\partial_2$, where $(C_*,\partial)$ is the cellular
complex (with respect to any \CW-structure) computing the homology of
$\FC\Ga$.
Indeed, the quotient $(C_1/\Im\partial_2)/H_1(\FC\Ga)$ is a subgroup of the
free abelian group $C_0$.

\subsection{The Alexander module}\label{s.AM}
Given a topological space~$X$ and an epimorphism $\Ga\:\pi_1(X)\onto G$
onto an abelian group~$G$, the homology of the covering $\tX\to X$
defined by~$\Ga$ are naturally $\Z[G]$-modules, $G$ acting by the deck
translations of the covering. The $\Z[G]$-module $H_1(\tX)$ is called the
\emph{Alexander module} of~$X$ or, more precisely, of pair $(X,\Ga)$. The
Alexander module depends only on
the group $\pi_1(X)$ and epimorphism~$\Ga$; algebraically, it is
the abelian group $\Ker\Ga/[\Ker\Ga,\Ker\Ga]$ on which
$G=\pi_1(X)/\Ker\Ga$ acts by
conjugation.

We employ the concept of Alexander module to compute $H_1(\FC\Ga)$.
First, let $\Ga$ be the identity map $0\:\GG\onto\GG/0\GG$
(awkward as it seems, this notation agrees with our convention;
we also have $1\:\GG\onto\GG/\GG=\{1\}$)
and consider the
ring $\GL:=\Z[\GG]$. Note that, unlike~$\FF0$, the unramified
covering $\FC0$ still
makes sense.
The group $\pi_1(\FC1)$ is computed using Zariski--van Kampen
theorem~\cite{vanKampen} (this computation is essentially shown
in \autoref{fig.divisor}, see \cite{degt:Fermat} for the relations
and further details), and the map
$\pi_1(\FC1)\onto\GG$ is $\hg_1\mapsto t_1$, $\vg_2\mapsto t_2$,
$\vg_3\mapsto t_3$, $\vg_1,\vg_4,\hg_2\mapsto 1$.
Then, the $\GL$-module $H_1(\FC0)$ is found by means of
the Fox free calculus~\cite{Fox:I}.
It is more convenient to work with the module
\[*
\Am0:=C_1[0]/\Im\partial_2,
\]
where $(C_*[0],\partial)$ is an appropriate cellular
complex computing the homology
(or even just the fundamental group)
of $\FC0$; as explained in
\autoref{s.torsion}, that would suffice for our purposes.
As a $\GL$-module, $\Am0$ is generated by six elements
$\ag_1$, $\ag_2$, $\ag_3$, $\cg_1$, $\cg_2$, $\cg_3$
(corresponding, in the order listed, to the six generators
$\hg_1$, $\vg_2$, $\vg_3$, $\hg_2$, $\vg_4$, $\vg_1$ of $\pi_1(\FC1)$ shown in
\autoref{fig.divisor}),
which are subject to the six relations
\begin{gather*}
(t_2t_3-1)\cg_1=(t_1t_3-1)\cg_2=(t_1t_2-1)\cg_3=0,\label{rel.bb}\\
 (t_3-1)\cg_1+(t_3-1)\ag_2-(t_2-1)\ag_3=0,\label{rel.a1}\\
 (t_3-1)\cg_2+(t_3-1)\ag_1-(t_1-1)\ag_3=0,\label{rel.a2}\\
 (t_1-1)\cg_3+(t_1-1)\ag_2-(t_2-1)\ag_1=0.\label{rel.a3}
\end{gather*}
Recall that we also have $t_0t_1t_2t_3=1$ in~$\GL$.

Now, given a finite quotient $\Ga\:\GG\to G$, the group ring
$\GL[\Ga]:=\Z[G]$ is naturally a $\GL$-module and the complex
$(C_*[\Ga],\partial)$ is obviously $(C_*[0],\partial)\otimes_\GL\GL[\Ga]$.
Hence, the new
$\GL[\Ga]$-module $\Am\Ga$ is $\Am0\otimes_\GL\GL[\Ga]$. In some cases, it is
more convenient to work with the submodule $\Bm\Ga\subset\Am\Ga$
generated by $\cg_1,\cg_2,\cg_3$.
(Note though that it is not always easy to find the defining relations for
$\Bm\Ga$.)
Since the complex
\[*
0\longto\Am\Ga/\Bm\Ga\longto\GL[\Ga]\longto0
\]
computes the homology $H_0=\Z$ and
$H_1=\Ker\Ga\subset\GG$ of the space $\FF\Ga\sminus\pi_A\1R$,
the two modules
have the same integral torsion. Summarizing, we arrive at the following
algebraic description of $\KF\Ga$ and $\TF\Ga$.

\theorem[see~\cite{degt:Delsarte}]\label{th.Delsarte}
For any finite quotient $\Ga\:\GG\to G$ one has
\[*
\aligned
\TF\Ga&=\Ext_\Z(\Am\Ga,\Z)=\Ext_\Z(\Bm\Ga,\Z),\\
\rank_\Z\KF\Ga&=\rank_\Z\Am\Ga-\ls|G|+1=\rank_\Z\Bm\Ga+3.
\endaligned
\]
\endtheorem

In other words, the torsion $\TF\Ga$ in question
is isomorphic to the integral torsion of either of
the two modules
$\Am\Ga$ or $\Bm\Ga$.

Note also that, even if $\TF\Ga\ne0$, a sufficiently good description of
this group would still let one recover the complete structure of
$\NS(\FF\Ga)$. For example, one can use the technique of discriminant forms
introduced in~\cite{Nikulin:forms}.
From this point of view,
\autoref{th.Delsarte} does give us a suitable description of $\TF\Ga$, as it
actually
places this group to the discriminant group $\SF\Ga\dual\!/\SF\Ga$.

\subsection{Vanishing and bounds}\label{s.partial}
Numeric experiments with random matrices show that, typically, $\TF\Ga\ne0$,
even if $\gcd(m,6)=1$ (see \autoref{s.Delsarte}).
However, the vanishing of the group $\TF\Ga$ can be established
in several important special cases. We have the following theorem.

\theorem[see \cite{degt:Delsarte,degt:Fermat}]\label{th.=0}
In each of the following three cases, one has $\TF\Ga=0$\rom:
\roster
\item\label{=0.Fermat}
$\FF\Ga$ is a Fermat surface, \ie, $\Ga=m\in\N_+$\rom;
\item\label{=0.cyclic}
$\FF\Ga$ is \emph{cyclic}, \ie, the image~$G$ of~$\Ga$ is a cyclic
group\rom;
\item\label{=0.unramified}
$\FF\Ga$ is \emph{unramified at~$\infty$}, \ie, $\Ga(t_0)=1$.
\endroster
\endtheorem

Statement~\iref{=0.unramified} in \autoref{th.=0} was a toy example
considered in~\cite{degt:Fermat}. Statement~\iref{=0.cyclic} is proved
in~\cite{degt:Delsarte} by comparing the dimensions
$\dim_\Bbbk\Am\Ga\otimes\Bbbk$, where $\Bbbk$ is either~$\C$ or a finite
field $\Bbb{F}_p$: if $G$ is cyclic, $\GL[\Ga]\otimes\Bbbk=\Bbbk[G]$
is a principal
ideal domain and the dimension of a module can be computed
algorithmically using elementary
divisors of the matrix of relations.

Statement~\iref{=0.Fermat} is more involved. In~\cite{degt:Fermat}, it is
proved by considering an appropriate rather long filtration
\[*
0=A_0\subset A_1\subset\ldots\subset A_7=\Am\Ga
\]
and \emph{estimating}
from above the
length $\ell(A_{i+1}/A_i)$ of each quotient.
(Recall that the \emph{length}
$\ell(A)$ of an abelian group~$A$ is the minimal number of generators
of~$A$, whereas its \emph{rank} $\rank A$ is the maximal number of linearly
independent elements.
Always $\rank A\le\ell(A)$, and a
finitely generated abelian group~$A$ is free if and only if
$\ell(A)=\rank A$.) Luckily, these estimates sum up to the expected rank
$\rank\Am\Ga$ given by~\eqref{eq.rank.S} and \autoref{th.Delsarte}; hence,
each quotient $A_{i+1}/A_i$ is a free abelian group,
and so is the original module $\Am\Ga$.

The same approach can be used in the general case, but the counts no longer
match; hance, we only obtain an estimate on the size of $\TF\Ga$.
To state the next theorem, we need to introduce a few invariants measuring
the non-uniformity of the homomorphism~$\Ga$.
(Note that the group~$\GG$ is to be considered in its canonical generating
set $t_0,t_1,t_2,t_3$
introduced in \autoref{s.Delsarte}; the only automorphisms allowed are
permutations of the generators. This rigidity explains also why we
are using four generators instead of three.)
First, consider the following subgroups
of~$\GG$:
\roster*
\item
$\GG_{ij}$, generated by $t_i$ and~$t_j$, $i,j=0,1,2,3$;
\item
$\GG_i$, generated by $t_it_j$ and $t_it_k$, $i=1,2,3$ and
$\{i,j,k\}=\{1,2,3\}$;
\item
$\GG\two:=\sum_i\GG_i$, generated by $t_1t_2$, $t_1t_3$, and $t_2t_3$.
\endroster
In more symmetric terms,
$\GG_i$ depends only on the partition $\{\{0,i\},\{j,k\}\}$ of the index
set, see~\eqref{eq.partition},
and $\GG\two$ is generated by all products $t_it_j$, $i,j=0,1,2,3$;
one has $[\GG:\GG\two]=2$.
Now, for a finite quotient $\Ga\:\GG\onto G$, denote $G_*:=G/\Ga(\GG_*)$
(where the subscript $*$ is one of the symbols
$ij$, $i$, or~$=$ as above) and define
$\Gd\of\Ga:=\ls|G\two|-1\in\{0,1\}$.
Let, further, $\exp G$ be the minimal positive integer~$m$ such that
$mG=0$. (This notion applies to any abelian group~$G$; in our case, it is
also the minimal positive integer~$m$ such that $m\GG\subset\Ker\Ga$).

In this notation, the fundamental group
$\pi_1(\FF\Ga)$ found in~\cite{degt:Fermat} is given by
\[
\pi_1(\FF\Ga)=H_1(\FF\Ga)=\Ker\Ga\big/\prod(\GG_{ij}\cap\Ker\Ga),
\label{eq.pi1}
\]
the product running over all pairs $0\le i<j\le3$.
This group is trivial in any of the three special classes considered in
\autoref{th.=0}.
In general, as shown in~\cite{degt:Delsarte},
the group $\pi_1(\FF\Ga)$ is cyclic
and its order $\ls|\pi_1(\FF\Ga)|$ divides the \emph{height}
$\height\Ga:=\exp G/n$, where $n$ is the maximal integer such that
$\Ker\Ga\subset n\GG$.

\theorem[see \cite{degt:Delsarte}]\label{th.bounds}
For any finite quotient $\Ga\:\GG\onto G$,
one has
\[*
\rank\KF\Ga=
 \sum_{0\le i<j\le3}\ls|G_{ij}|+\sum_{1\le i\le3}\ls|G_{i}|
 -3-\Gd\of\Ga.
\]
Besides,
$\ell(\TF\Ga)\le6+\Gd\of\Ga$ and
$\exp\TF\Ga$ divides $(\exp G)^3\!/\ls|G|$.
\endtheorem

The bound on $\ell(\TF\Ga)$ is sharp, whereas that on $\exp\TF\Ga$ is
probably not.

Analyzing the proof, one can also establish the almost vanishing of the
torsion in the case of \emph{Brieskorn surfaces} (called \emph{diagonal
Delsarte surfaces} in~\cite{degt:Delsarte}), \ie, those given by an affine
equation of the form
\[*
x^{m_1}+y^{m_2}+z^{m_3}=1,
\]
so that $\Ga$ is the projection
$\GG\onto\GG/(t_1^{m_1}=t_2^{m_2}=t_3^{m_3}=1)$.
For such surfaces, $\TF\Ga$ is cyclic: one has
$\ell(\TF\Ga)\le\Gd\of\Ga$
and the order $\ls|\TF\Ga|$ divides the ratio
\[*
h(m_1,m_2,m_3):=
\frac{\lcm_{1\le i<j\le3}\bigl(\gcd(m_i,m_j)\bigr)}{\gcd(m_1,m_2,m_3)}
=\sqrt{\frac{\mathstrut{m_1m_2m_3}}{\gcd(m_1,m_2,m_3)^{3}}}.
\]
(Observe that $h(m_1,m_2,m_3)=1$ if $m_1=m_2=m_3$, \ie, in the case
of classical Fermat surfaces.)
Examples show that these bounds are also sharp.

\section{Higher dimensions}

\subsection{The reduction}\label{s.reduction}
Now, let us consider a Fermat variety $\Fermat:=\Fermat_m^d$ of even
dimension $d=2k\ge4$. Denote by $\CJ:=\CJ(d)$ the set of partitions as
in~\eqref{eq.partition}; each element $J\in\CJ$ gives rise to $m^{k+1}$
subspaces in~$\Fermat$. To put the statements in a slightly more general
form, we will pick a nonempty subset $\CK\subset\CJ$ and denote by
$\dv_\CK\subset\Fermat$ the union of all subspaces $L_{J,\eta}$,
see~\eqref{eq.L}, with $J\in\CK$ and $\eta$ running over all sequences of
roots of~$(-1)$.
Denoting by $\inj\:\dv_\CK\into\Fermat$ the inclusion, consider the
groups
\[*
\sset_\CK:=\Im\bigl[\inj_*\:H_d(\dv_\CK)\to H_d(\Fermat)\bigr],\qquad
\tset_\CK:=\Tors\bigl(H_d(\Fermat)/\sset_\CK\bigr).
\]
Clearly, $\sset_\CJ=\sset_m^d$ and $\tset_\CJ=\tset_m^d$.

In what follows, we always regard $H_d(\Fermat)$ as a unimodular lattice by
means of the Poincar\'e duality isomorphism $H_d(\Fermat)\to
H^d(\Fermat)=H_d(\Fermat)\dual$. (Here and below, we denote by
$A\dual:=\Hom(A,\Z)$ the dual of an abelian group~$A$.)
Consider the subspace $Y:=\Fermat\sminus\{z_0=0\}$ and
let $\Hc\subset H_d(\Fermat)$ be the image of the inclusion homomorphism
$H_d(Y)\to H_d(\Fermat)$; it coincides with the
orthogonal complement of the class $h\in H_d(\Fermat)$ of the intersection of
$\Fermat$ with a generic $(d+1)$-plane. Since the lattice $H_d(\Fermat)$ is
unimodular, the composition $H_d(\Fermat)=H_d(\Fermat)\dual\onto\Hc\dual$ is
surjective; in fact, we have a short exact sequence
\[*
0\longto\Z h\longto H_d(\Fermat)\longto\Hc\dual\longto0.
\]
Let $\sset_\CK'\subset\Hc\dual$ be
the image of~$\sset_\CK$. Since
$h\in\sset_\CK$, we have $H_d(\Fermat)/\sset_\CK=\Hc\dual\!/\sset_\CK'$
and $\rank\sset_\CK'=\rank\sset_\CK-1$.
(Recall that we assume that $\CK\ne\varnothing$; hence, fixing $J\in\CK$ and
all but one $\eta_i$ in~\eqref{eq.L}, we obtain $m$ spaces whose classes sum
up to~$h$.)

Let $\GL:=\Z[G_m]$, see~\eqref{eq.Gm}. To make the notation more
conventional, we rename the canonical generators
$(1,\ldots,\exp(2\pi i/m),\ldots,1)$
of~$G_m$ into
$t_0,\ldots,t_{d+1}$ and regard~$\GL$ as the quotient of the ring
$\Z[t_0^{\pm1},\ldots,t_{d+1}^{\pm1}]$ of Laurent polynomials by the ideal
generated by $t_0\ldots t_{d+1}-1$ and $t_i^m-1$, $i=0,\ldots,d+1$.
Since $G_m$ acts on~$\Fermat$, all homology groups involved are naturally
$\GL$-modules and,
$G_m$ acting identically on the fundamental class $[\Fermat]$,
the lattice structure on $H_d(\Fermat)$ is $\GL$-invariant.

For further statements, we need to prepare several polynomials. Let
\[*
\Gf(t):=\sum_{i=0}^{m-1}t^i=\frac{t^m-1}{t-1},\qquad
\Gr(x,y):=\sum_{0\le i\le j\le m-2}x^jy^i
\]
and, for $J\in\CJ$ as in~\eqref{eq.partition}, denote
\[*
\tau_J:=\prod_{i=0}^k(t_{q_i}-1),\quad
\psi_J:=\tau_J\prod_{i=1}^k\Gf(t_{p_i}t_{q_i}),\quad
\rho_J:=\prod_{i=1}^k\rho(t_{p_i},t_{q_i}).
\]
Also, for any quotient ring~$R$ of~$\GL$, including $\GL$ itself,
we will denote by $\bar R$ its
``reduced'' version, \viz.
$\bar R:=R\big/\sum_{i=0}^{d+1}R\Gf(t_i)$.

The advantage of using~$Y$ instead of~$\Fermat$ is the fact that this space
has extremely simple homology, which have been extensively studied as the
vanishing cycles of a Pham--Brieskorn singularity.
Fix a number $\zeta\in\C$ such that $\zeta^m=-1$ and consider the
topological simplex
\[*
\Delta:=\bigl\{(s_1,\ldots,s_{d+1})\zeta\bigm|
 s_s\in[0,1],\ s_1^m+\ldots+s_{d+1}^m=1\bigr\}\subset Y.
\]
Then, the so-called \emph{Pham polyhedron}
\[*
\Sigma:=(1-t_1\1)\ldots(1-t_{d+1}\1)\Delta
\]
is a cycle in~$Y$; in fact, $\Sigma$ is a topological sphere.

\theorem[see~\cite{Pham}]
The group $H_d(Y)$ is the free $\bar\GL$-module generated by~$\Sigma$.
\endtheorem

Therefore, $\Hc\dual$ is an ideal in $\GL=\GL\dual$
(where all groups dual to $\GL$-modules are regarded as
$\GL$-modules with respect to the
contragredient $G$-action)
and, in order to find the
image $\sset_\CK'\subset\Hc\dual\subset\GL$,
it suffices to compute the algebraic intersection of
each space~$L_{J,\eta}$ with~$\Sigma$. This is done in~\cite{degt.Shimada},
and, omitting intermediate details, the
result can be stated as follows.

\theorem[see~\cite{degt.Shimada}]\label{th.Shimada}
For each $J\in\CJ$, one has $\sset_{J}'=\GL\psi_J\subset\GL$. Hence, for a
subset
$\CK\subset\CJ$, one has $\sset_\CK'=\sum_J\GL\psi_J\subset\GL$, the
summation running over $J\in\CK$.
\endtheorem

\subsection{Partial vanishing statements}\label{s.vanishing}
Using \autoref{th.Shimada} and employing various dualities and torsion-free
quotients, we can obtain several expressions for~$\tset_\CK$. In the
statement below, for $J\in\CJ$ as in~\eqref{eq.partition}, we use the ring
\[*
\textstyle
\GL_J:=\GL\big/\sum_i\GL(t_{p_i}t_{q_i}-1),\quad i=0,\ldots,k+1,
\]
and $1_J$ stands for the unit in $\GL_J$ or $\bar\GL_J$.

\theorem[see~\cite{degt.Shimada}]\label{th.TK}
Let $\CK\subset\CJ$ be a nonempty subset. Then, the torsion $\tset_\CK$ is
isomorphic to the integral torsion of any of the following modules\rom:
\roster
\item\label{T.direct}
the ring $\GL\big/\sum_J\GL\psi_J$, $J\in\CK$\rom;
\item\label{barT.direct}
the ring $\bar\GL\big/\sum_J\bar\GL\rho_J$, $J\in\CK$\rom;
\item\label{T.dual}
the $\GL$-module $M_\CK:=\bigl(\bigoplus_J\GL_J\bigr)/\GL\tau$,
where $\tau:=\sum_J\tau_J1_J$ and $J\in\CK$\rom;
\item\label{barT.dual}
the $\bar\GL$-module $\bar M_\CK:=\bigl(\bigoplus_J\bar\GL_J\bigr)/\bar\GL1$,
where $1:=\sum_J1_J$ and $J\in\CK$.
\endroster
\endtheorem

Denoting by~$T$ the integral torsion of the respective module in
\autoref{th.TK},
in Statements~\iref{T.direct} and~\iref{barT.direct} we have canonical
isomorphisms $\tset_\CK=T$, whereas in~\iref{T.dual}
and~\iref{barT.dual}, the isomorphisms are $\tset_\CK=\Hom(T,\Q/\Z)$.

If $d=2$, the module $M_\CJ$ in \autoref{th.TK}\iref{T.dual} coincides with
the module~$\Bm{m}$
introduced in \autoref{s.AM}.
Both $M_\CJ$ and $\bar M_\CJ$ appear as intermediate quotients of
the filtration used in the proof of \autoref{th.=0}\iref{=0.Fermat}.

\conjecture\label{conj.TK}
For the full set $\CK=\CJ$, one has $\tset_\CJ=0$.
\endconjecture

This conjecture holds true in small dimensions $d=0$ (obvious) and $d=2$
(see \autoref{th.=0}) and is supported by some numerical evidence: by a
computer aided computation, we managed to establish the vanishing of
$\tset_\CJ$ for the values
\[*
(d,m)=(4,m),\ 3\le m\le12,\quad
(6,3),\ (6,4),\ (6,5),\ \text{and}\ (8,3).
\]
Unfortunately, we failed to prove the conjecture in full generality. It is
not difficult to show (see~\cite{degt.Shimada}) that
$\gcd(\ls|\tset_\CK|,p)=1$ for each prime $p\nmid m$. One can also show that
$\tset_\CK=0$ for some special subsets $\CK\subset\CJ$.
As an important example (which can probably be used as a base for induction),
fixing the degree~$m$ and denoting by $(\cdot)$ the dependence on the dimension,
consider the natural inclusions $\CJ(s)\subset\CJ(d)$, $s=2l\le d=2k$,
extending each partition $J\in\CJ(s)$ \emph{identically beyond~$s$}, \ie,
attaching the pairs $\{2i,2i+1\}$, $i=l+1,\ldots,k$.
Then, given $\CK\subset\CJ(s)$, for the module~$\bar M_\CK$ in
\autoref{th.TK}\iref{barT.dual} one can easily see that
\[*
\bar M_\CK(d)=\bar M_\CK(s)\otimes_\Z
 \bar\Delta_{s+2}\otimes_\Z\bar\Delta_{s+4}\otimes_\Z\ldots\otimes_\Z\bar\Delta_d,
\]
where $\bar\Delta_i:=\Z[t_i^{\pm1}]/\phi(t_i)$. Hence, we have
stabilization
\[
\tset_\CK(d)=\tset_\CK(s)\otimes_\Z
 \bar\Delta_{s+2}\otimes_\Z\bar\Delta_{s+4}\otimes_\Z\ldots\otimes_\Z\bar\Delta_d.
\label{eq.stab}
\]
In particular, $\tset_\CK(d)=0$ if and only if $\tset_\CK(s)=0$.
This observation applies to the image $\CJ(s,d)$ of $\CJ(s)$ in $\CJ(d)$;
thus
$\tset_{\CJ(s,d)}=0$ for all dimensions $d\ge s$
if and only if $\tset_{\CJ(s)}=0$. As a consequence, for
any~$d$,
\[
\tset_{\CJ(0,d)}=\tset_{\CJ(2,d)}=0;
\label{eq.0,2}
\]
the former is obvious, and the latter follows from
\autoref{th.=0}\iref{=0.Fermat}.

\subsection{Other classes of varieties}\label{s.others}
The last two statements~\eqref{eq.stab}, \eqref{eq.0,2} have
a geometric interpretation.
Given $s=2l<d=2k$, consider the
\emph{partial Fermat variety} $\Fermat_m^{s,d}$ given by an equation
of the form
\[*
f_0(z_0,z_1)+f_1(z_2,z_3)+\ldots+f_k(z_d,z_{d+1})=0,
\]
where each $f_i$ is a homogeneous bivariate polynomial of degree~$m$ and
\[*
f_0(u,v)=f_1(u,v)=\ldots=f_l(u,v)=u^m+v^m,
\]
whereas all other terms
are generic, pairwise distinct, and other than $u^m+v^m$. Arguing as in
\autoref{s.planes}, one can see that $\Fermat_m^{s,d}$ contains several group
of $d$-spaces: each group consists of $m^{k+1}$ spaces, and the groups are
indexed by the members of $\CJ(s,d)$. Now, observe that the proof of
\autoref{th.TK} is purely topological; hence, we can deform $\Fermat_m^{s,d}$
to~$\Fermat_m^d$ (followed by a deformation of the $d$-spaces
in $\Fermat_m^{s,d}$ to some of those
in~$\Fermat_m^d$) and apply \autoref{th.TK}, obtaining the following corollary.

\corollary[see~\cite{degt.Shimada}]\label{cor.partial}
The subspaces contained in $\Fermat_m^{s,d}$ generate a primitive
subgroup in the N\'eron--Severi lattice $\NS(\Fermat_m^{s,d})$ for any
dimension
$d\ge s$ if and only if they do so for $d=s$, \ie, if $\tset_m^s=0$.
\endcorollary

\corollary[see~\cite{degt.Shimada}]\label{cor.0,2}
For $s=0$
or~$2$, the subspaces contained in $\Fermat_m^{s,d}$ generate a primitive
subgroup in the N\'eron--Severi lattice $\NS(\Fermat_m^{s,d})$.
\endcorollary

If $s=0$, we can also choose $f_0$ generic, retaining the $m^{k+1}$ spaces
contained in $\Fermat_m^{0,d}$.
In this case, if $d=2$ (lines in surfaces) and $m$ is
prime, the $m^2$ lines contained in $\Fermat_m^{0,2}$ are known to generate
the rational N\'eron--Severi lattice
$\NS(\Fermat_m^{0,2};\Q)$, see, \eg,~\cite{Boissiere.Sarti}.
\autoref{cor.0,2} implies
that these lines generate $\NS(\Fermat_m^{0,2})$ over~$\Z$
(see~\cite{degt:Fermat}).

\section{Open problems}\label{S.problems}

Apart from \autoref{conj.TK}, there are a few other interesting open questions
that may be worth stating explicitly.

As explained in \autoref{s.partial}, typically, for a Delsarte surface
$\FF\Ga$ one has $\TF\Ga\ne0$.
Naturally, one may ask if there are other classes of surfaces for which one can
assert that $\TF\Ga=0$ or obtain a bound on the size of this group better
than that given by \autoref{th.bounds}.
In
\autoref{th.=0}, the
Delsarte surfaces are treated according to the complexity
(or non-uniformity)
of the finite quotient $\Ga\:\GG\onto G$.
However, there are other taxonomies which, from many points of view, may seem
much more natural. For example, one can classify Delsarte surfaces according to
the singularities of the original (not yet resolved) projective hypersurface
given by~\eqref{eq.Delsarte}.
Thus, it is known that
there are ten families (one of them being Fermat) of nonsingular Delsarte
surfaces, see~\cite{Kogure}, and
$83$ families of those with $\bA$--$\bD$--$\bE$ singularities,
see~\cite{Heijne}.
The Picard ranks for these families were computed in \cite{Heijne,Kogure}.

\problem
Does the vanishing $\TF\Ga=0$ hold for all nonsingular Delsarte surfaces
$\FF\Ga$?
For those with $\bA$--$\bD$--$\bE$ singularities?
\endproblem

\problem
Are there sharper bounds on the size
(length, order, exponent) of the group $\TF\Ga$ in terms of the singularities
of~$\FF\Ga$?
\endproblem

As another generalization, one can consider a Fermat surface $\Fermat_m^2$
of a degree~$m$ not prime to~$6$, so that the lines do \emph{not} generate
$\NS(\Fermat_m^2;\Q)$. In some cases, there are explicit
lists of additional generators.
Thus, found in~\cite{Shioda:Fermat1}, there is a list of
relatively simple curves, lying in quadrics, cubics, and quartics, that
compensate for the terms $24(m/3)^*$ and $48(m/2)^*$ in~\eqref{eq.dim.NS}.
As in \autoref{s.planes}, the generating property is established by comparing
the ranks; hence, the question whether these curves (together with the lines)
generate the integral group $\NS(\Fermat_m^2)$ remains open.

\problem[T.~Shioda]
For the known explicit generating sets $\CS$ of
the group $\NS(\Fermat_m^2;\Q)$, is it true that the curves
constituting~$\CS$ also generate $\NS(\Fermat_m^2)$ over the integers?
In other words, is it true that the subgroup
$\SS{\CS}:=\sum\Z[C]$, $C\in\CS$, is
primitive in $H_2(\Fermat_m^2)$?
If not, what is the torsion of $H_2(\Fermat_m^2)/\SS{\CS}$?
\endproblem

\let\.\DOTaccent
\def\cprime{$'$}
\bibliographystyle{amsplain}
\bibliography{degt}


\end{document}